\RequirePackage{fix-cm} 
\documentclass[a4paper,twoside,11pt,reqno]{amsart}

\usepackage[english]{babel}
\usepackage[latin1]{inputenc}

\usepackage{indentfirst,a4wide,verbatim}

\usepackage{amsmath,amsfonts,amssymb,amsgen,amsbsy,eucal,mathrsfs,dsfont,stmaryrd}

\usepackage[square,numbers]{natbib}
\usepackage{mathtools}

\usepackage[pagewise]{lineno}

\usepackage{hyperref}
\hypersetup{
colorlinks=true, 
breaklinks=true, 
urlcolor= blue, 
linkcolor= blue, 
bookmarksopen=true, 
pdftitle={}, 
pdfauthor={}, 
pdfsubject={} 
}

\newcommand{\ud}{\mathrm{d}}

\newcommand{\half}{{\textstyle{1\over2}}}

\newtheorem{thm}{Theorem}
\newtheorem{lem}{Lemma}

\usepackage{enumitem}
\newlist{steps}{enumerate}{1}
\setlist[steps, 1]{label = Step \arabic*:}

\newcommand{\eqdef}{\stackrel{\text{\tiny{def}}}{=}}

\title[]{Local well-posedness of a Hamiltonian regularisation of the Saint-Venant system with uneven bottom}

\author[GUELMAME et al.]{Billel Guelmame, Didier Clamond and St\'ephane Junca}

\newcommand{\nfont}{\fontshape{n}\selectfont}

\address{({\nfont\textbf{Billel Guelmame}})  LJAD,  Inria \& CNRS,  Universit\'e C\^ote d'Azur, France.} 
\email{billel.guelmame@univ-cotedazur.fr}

\address{({\nfont\textbf{Didier Clamond}})  LJAD,  CNRS,  Universit\'e C\^ote d'Azur,  France.} 
\email{didier.clamond@univ-cotedazur.fr}

\address{({\nfont\textbf{St\'ephane Junca}})  LJAD,  Inria \& CNRS,   Universit\'e C\^ote d'Azur, France.} 
\email{stephane.junca@univ-cotedazur.fr}

\let\oldtocsection=\tocsection
 
\let\oldtocsubsection=\tocsubsection

\renewcommand{\tocsection}[2]{\hspace{0em}\oldtocsection {#1}{#2}}
\renewcommand{\tocsubsection}[2]{\hspace{2em}\oldtocsubsection{#1}{#2}}
\begin{document}

\begin{abstract}
We prove in this note the local (in time) well-posedness of a broad class of $2 \times 2$ symmetrisable 
hyperbolic system involving additional non-local terms.
The latest result implies the local well-posedness of the non dispersive regularisation of the 
Saint-Venant system with uneven bottom introduced by \citet{clamond2019hamiltonian}.
We also prove that, as long as the first derivatives are bounded, singularities cannot appear.
\end{abstract}

\maketitle


\medskip



\section{Introduction and main results}
\citet{ClamondDutykh2018a} have recently proposed  a new Hamiltonian regularisation of the Saint-Venant 
(rSV) system with a constant bottom. 
Such a regularisation avoids some defects of diffusive approximations which flatten the shocks, and the dispersive ones which add spurious oscillations.
This Hamiltonian regularisation   \cite{ClamondDutykh2018a} is  of a new kind, non-diffusive and non-dispersive. 
This new model has been mathematically studied in 
\citep{liu2019well,PuEtAl2018}.  
Inspired by \cite{ClamondDutykh2018a}, similar regularisations have been proposed for the inviscid Burgers 
equation \citep{guelmame2020global}, the scalar conservation laws \citep{guelmame2020hamiltonian}
 and the barotropic Euler system \citep{guelmame2020Euler}. 
A regularisation of the Saint-Venant equations with uneven bottom (rSVub) has also been proposed by 
\citet{clamond2019hamiltonian}. The latter equations, for the conservation of mass and momentum, can 
be written in the conservative form
\begin{subequations}\label{rSVmb}
\begin{align}
h_t\ +\,\left[\,h\,u\,\right]_x\ &=\  0, \\
\left[\,h\,u\,\right]_t\ +\,\left[\,h\,u^2\,+\, \half\, g\, h^2\,
+\,\varepsilon\,\mathscr{R}\,\right]_x\ 
&=\  \varepsilon\, g\, h^2\, \eta_x\, d_{xx}\ +\ g\, h\, d_x,\\ \label{Rdef}
\mathscr{R}\ \eqdef\, 2\, h^3\, u_x^{\,2}\ -\ h^3\,
\left[\,u_t\,+\,u\,u_x\,+\, g\,\eta_x\,\right]_x\ &-\ \half\, g\, h^2\, \left(\eta_x^{\,2} +\, 
2\, \eta_x\, d_x \right).
\end{align}
\end{subequations}
Here, $u=u(t,x)$ is the depth-averaged horizontal velocity, $h= h(t,x) \eqdef \eta(t,x) + d(t,x)$ denotes 
the total water depth, $\eta$ being the surface elevation from rest and $d$ being the water depth for 
the unperturbed free surface. 
We can assume, without losing generality via a change of frame of reference, that the spacial average 
of the depth  $\bar{d}$ is constant in time. In that case, the gravity acceleration $g=g(t)$ may be a 
function of time. Introducing the Sturm--Liouville operator 
\begin{equation}
\mathcal{L}_h\ \eqdef\ h\ -\  \varepsilon\, \partial_x\, h^3\, \partial_x,
\end{equation}
if $h>0$, the operator $\mathcal{L}_h$ is invertible, then the system \eqref{rSVmb} can be written on the form
\begin{subequations}\label{rSVmb2}
\begin{flalign}\label{mc}
h_t\ +\,\left[\,h\,u\,\right]_x\ &=\ 0, &&\\ \label{momG3}
u_t\, +\, u\,u_x\, +\, g\, \eta_x\
&=\ \varepsilon\, g\, \mathcal{L}_h ^{-1} \left\{ h^2\, \eta_x\, d_{xx} \right\}
-\, \varepsilon\, \mathcal{L}_h ^{-1}\, \partial_x\,  \left\{ 2\, h^3\, u_x^{\,2}\,  
-\, \half\, g\, h^2\, \left(\eta_x^{\,2} +\, 2\, \eta_x\, d_x \right)\right\}. \hspace{-1cm} &&
\end{flalign}
\end{subequations}

The regularised Saint-Venant system admits a Hamiltonian structure, thus, it necessarily conserves\footnote{If the energy changes, it is necessarily due to exterior forces (the moving bottom) and not due to the dynamic itself.} the corresponding energy for smooth solutions.
The energy equation of \eqref{rSVmb2} writes 
\begin{gather}
\left[\,\half\,h\,u^2\,+\,\half\,\varepsilon\,h^3\,u_x^{\,2}\,+\, \half\, g\, \eta^2\,+\, \half\,
\varepsilon\, g\, h^2\, \eta_x^{\,2}\,\right]_t\ \nonumber\\
+\ \left[\/\left(\,\half\,h\,u^2\,+\, g\, h\,\eta\,+\, \half\, \varepsilon\, h^3\,u_x^{\,2}\,
+\, \half\,  \varepsilon\, g\, h^2\,\eta_x^{\,2}\, +\,\varepsilon\,\mathscr{R}\,\right)u\,
+\,\varepsilon\, g\, h^3\,\eta_x\,u_x\,\right]_x\  \nonumber\\
=\ \half\, \dot{g} \left( \eta^2\, +\, \varepsilon\, h^2\, \eta_x^2 \right)\ -\  g\, \eta\, d_t\ 
-\ \varepsilon\, g\, h^2\, \eta_x\, d_{xt}. \label{eneSV}
\end{gather}
Note that, injecting \eqref{momG3} in \eqref{Rdef}, one obtains the alternative definition of $\mathscr{R}$
\begin{equation*}
\mathscr{R}\ =\
\left( 1\, +\, \varepsilon\, h^3\, \partial_x\, \mathcal{L}_h^{-1}\, \partial_x \right)\left\{ 
2\,h^3\, u_x^{\,2}\,-\, \half\, g\,h^2\, \left(\eta_x^{\,2}\/+\/2\/\eta_x\/d_x \right) 
\right\}\,-\ \varepsilon\, h^3\, \partial_x\, \mathcal{L}_h^{-1}\, \left\{g\, h^2\, \eta_x\, 
d_{xx} \right\}.
\end{equation*}

The rSV and rSVub equations can be compared with the Serre--Green--Naghdi and the two-component 
Camassa--Holm equations. The local well-posedness of those equations have been studied in the literature 
(see, e.g., \citep{Israwi2011,guan2010well}; see also \citep{coclite2009well} for higher-order 
Camassa--Holm equations).
\citet{liu2019well} have proved the local well-posedness of the rSV equations introduced in 
\cite{ClamondDutykh2018a} for constant depth.  \citet{liu2019well}  have constructed some small 
initial data, such that the corresponding solutions blow-up in finite time.
The goal of the present note is to prove the local (in time) well-posedness of the rSVub equations. 
To this aim, we prove first the local (in time) well-posedness of a general $2 \times 2$ 
symmetrisable hyperbolic system. Then, using some estimates of the operator 
$\mathcal{L}_h^{-1}$, we prove that the system \eqref{rSVmb2} is locally well-posed in the Sobolev space $H^s(\mathds{R})$ 
for any real number $s \geqslant 2$.
We also prove that if the $L^\infty$-norm of the first derivatives remain bounded, then the 
singularities cannot appear in finite time.

In order to state the main results of this note, let ${d}$ be a smooth function of $t$ and $x$ with
\begin{equation}\label{ddef}
{h}\ \eqdef\ \eta\ +\ d, \qquad \bar{{d}}\ \eqdef\ \lim_{|x| \to +\infty} d(t,x)\ >\ 0, 
\qquad \mathrm{and} \qquad \inf_{(t,x) \in \mathds{R}^+ \times \mathds{R}} d(t,x)\ >\ 0,
\end{equation}
defining the Sobolev space $H^s \eqdef H^s(\mathds{R}) \eqdef \left\{u,\ \int_\mathds{R}\! \left( 1 + \xi^2 \right)^s \left| \hat{u}(\xi) \right|^2\, \mathrm{d} \xi < + \infty  \right\}$ where $\hat{u}$ is the Fourier transform of $u$, then
\begin{thm}\label{existenceSV}
Let $\tilde{m} \geqslant s \geqslant 2$, $0<g\in {C}^1([0,+\infty[)$, 
$d-\bar{d} \in {C}([0,+\infty], H^{s+1}) \cap {C}^1([0,+\infty], H^{s})$ and let $W_0=(\eta_0,u_0)^\top \in H^s$ satisfying $\inf_{x \in \mathds{R}} h_0(x) \geqslant h^* >0$, then there exist $T>0$ and a unique solution $W=(\eta,u)\in {C}([0,T], H^s) \cap {C}^1([0,T], H^{s-1})$ of \eqref{rSVmb2} satisfying the non-zero depth condition 
$ \inf_{(t,x) \in [0,T] \times \mathds{R}}\, h(t,x)\ >\ 0. $
Moreover, if the maximal time of existence $T_{max}< +\infty$, then 
\begin{equation}\label{buc1}
\lim_{t \to T_{\text{max}}} \|W\|_{H^s}\ =\ +\infty \qquad or \qquad \inf_{(t,x) \in 
[0,T_{\text{max}}[ \times \mathds{R}} h(t,x)\ =\ 0.
\end{equation}
\end{thm}
Using the energy equation \eqref{eneSV} and some estimates, the blow-up criteria \eqref{buc1} can be improved.
\begin{thm}\label{thm:blowup}
For any interval $[0,T] \subset [0,T_{\text{max}}[$ ($T_{\text{max}}$ is the life span of the smooth solution), 
there exists $C>0$, such that $\forall t\in [0,T]$ we have 
\begin{equation}\label{eb}
\mathscr{E}(t) \eqdef\ \int_\mathds{R} \left[\,\half\,h\,u^2\,+\,\half\,\varepsilon\,h^3\,u_x^{\,2}\,+\, 
\half\, g\, \eta^2\,+\, \half\,\varepsilon\, g\, h^2\, \eta_x^{\,2}\,\right]\, \ud x\ \leqslant\ C.
\end{equation}
Moreover, if $T_{\text{max}}<+\infty$, then
\begin{equation}
\lim_{t \to T_{max}} \|W_x\|_{L^\infty}\ =\ +\infty.
\end{equation}
\end{thm}

Section \ref{sec:gesys} is devoted to prove the local well-posedness of a general $2 \times 2$ system.
The proofs of Theorems \ref{existenceSV} and \ref{thm:blowup} are given in Section \ref{sec:proof}.

\section{Local well-posedness of a general $2 \times 2$ system}\label{sec:gesys}

We prove here the local well-posedness of a class of systems with non-local 
operators in the $H^s$ space with $s>3/2$. 
Let ${d}$ be a smooth function such that \eqref{ddef} holds. Let also $N \geqslant 1$ be a natural 
number and $G\eqdef(g_1,\cdots,g_N)$ be a smooth function of $t$ and $x$, possibly depending on $d$, 
such that 
\begin{equation}
g_\infty(t)\ \eqdef\ g_1(t,\infty)\ =\ \lim_{|x| \to \infty} g_1(t,x)\ >\ 0 \quad \text{and}  
\quad g_{\text{inf}}\ \eqdef\ \inf_{(t,x) \in \mathds{R}^+ \times \mathds{R}} g_1(t,x)\ >\ 0.
\end{equation}
Let $f(d,h)$ be a positive function and let $f_1,f_2$ be functions of $d$, $h$, $u$, 
$\eta_x$, $u_x$ and $G$. Let also $a,b,c,f_3,f_4$ be functions of $d$, $h$, $u$ and $G$. 
We consider the symmetrisable hyperbolic system
\begin{subequations}\label{system}
\begin{align}
\eta_t\ +\ a(d,{h},u,G)\, \eta_x\ +\ b(d,{h},u,G)\, u_x\ &=\ \mathfrak{A}_1\, f_1\ +\  
\mathfrak{A}_3\, f_3,  \\
u_t\ +\ g_1\, f(d,{h})\, b(d,{h},u,G)\, \eta_x\ +\ c(d,{h},u,G)\, u_x\ &=\ \mathfrak{A}_2\, f_2\ 
+\  \mathfrak{A}_4\, f_4,
\end{align}
\end{subequations}
where the $\mathfrak{A}_j$  are linear operators depending on $h$ and $u$. In order to obtain the 
well-posedness of the system \eqref{system} in $H^s$ with $s>3/2$, we define $W \eqdef (\eta,u)^T$, 
$G_0 \eqdef (g_\infty,0,\cdots,0)$ and
\begin{align*}
B(W)\  &\eqdef\ \begin{pmatrix} a(d,{h},u,G) & b(d,{h},u,G) \\ g_1\, f(d,{h})\, b(d,{h},u,G) &  c(d,{h},u,G) \end{pmatrix}, \\ 
F(W)\ &\eqdef\  \begin{pmatrix} \mathfrak{A}_1\, f_1(d,{h},u,{h}_x,u_x,G)\ +\  \mathfrak{A}_3\, f_3(d,{h},u,G) \\  \mathfrak{A}_2\, f_2(d,{h},u,{h}_x,u_x,G)\ +\  \mathfrak{A}_4\, f_4(d,{h},u,G) \end{pmatrix},
\end{align*}
the system \eqref{system} can be written as
\begin{equation}\label{generalsystem}
W_t\ +\ B(W)\, W_x\ =\ F(W), \qquad W(0,x)\ =\ W_0(x).
\end{equation} 
We assume that:
\begin{itemize}
\item[(A1)] For $s \leqslant \tilde{m} \in \mathds{N}$, we have 
\begin{itemize}
\item[$\bullet$] $d-\bar{{d}},g_1-g_\infty,g_2,g_3,\cdots,g_N \in {C}(\mathds{R}^+,H^s)$ and $d-\bar{{d}},g_1-g_\infty 
\in {C}^1(\mathds{R}^+,H^{s-1})$;
\item[$\bullet$] $f \in {C}^{\tilde{m}+2}(]0,+\infty[^2)$ and for all ${h}_1,{h}_2>0$ we have $f({h}_1,{h}_2)>0$;
\item[$\bullet$] $f_1,f_2 \in {C}^{\tilde{m}+2}(]0,+\infty[^2 \times \mathds{R}^3 \times ]0,+\infty[ \times 
\mathds{R}^{N-1} )$;
\item[$\bullet$] $a,b,c,f_3,f_4 \in {C}^{\tilde{m}+2}(]0,+\infty[^2 \times \mathds{R} \times ]0,+\infty[ \times 
\mathds{R}^{N-1} )$;
\item[$\bullet$] $f_1(\bar{{d}},\bar{{d}},0,0,0,G_0)=f_2(\bar{{d}},\bar{{d}},0,0,0,G_0)=f_3(\bar{{d}},
\bar{{d}},0,G_0)=f_4(\bar{{d}},\bar{{d}},0,G_0)=0$.
\end{itemize}
\item[(A2)] For all $r \in [s-1,s]$, if $\phi \in H^r$ and $\psi \in H^{r-1}$, then 
\begin{align*}
\|\mathfrak{A}_1\, \psi \|_{H^r}\ +\ \|\mathfrak{A}_2\, \psi \|_{H^r}\ &\leqslant\ C(s,r,d,\|W\|_{H^r})\, 
\|\psi\|_{H^{r-1}},\\
\|\mathfrak{A}_3\, \phi \|_{H^r}\ +\ \|\mathfrak{A}_4\, \phi \|_{H^r}\ &\leqslant\ C(s,r,d,\|W\|_{H^r})\, 
\|\phi\|_{H^{r}}.
\end{align*}
\item[(A3)] 
If $\phi,W,\tilde{W} \in H^s$ and $\psi \in H^{s-1}$, then
\begin{align*}
\| (\mathfrak{A}_1(W)- \mathfrak{A}_1(\tilde{W}))\, \psi \|_{H^{s-1}}\ +\ \| (\mathfrak{A}_2(W)- 
\mathfrak{A}_2(\tilde{W}))\, \psi \|_{H^{s-1}}\ &\leqslant\ C\, \|W- \tilde{W}\|_{H^{s-1}},\\
\| (\mathfrak{A}_3(W)- \mathfrak{A}_3(\tilde{W}))\, \phi \|_{H^{s-1}}\ +\ \| (\mathfrak{A}_4(W)- 
\mathfrak{A}_4(\tilde{W}))\, \phi \|_{H^{s-1}}\ &\leqslant\ C\, \|W- \tilde{W}\|_{H^{s-1}},
\end{align*}
where $C=C\left(s,d,\|W\|_{H^s},\|\tilde{W}\|_{H^s},\|\phi\|_{H^s},\|\psi\|_{H^{s-1}} \right)$.
\end{itemize}
Note that if ${h}$ is far from zero (i.e., $\inf h >0$), then $g_1 f(d,{h})$ is positive and 
far from zero. Then, the system \eqref{system} is symmetrisable and hyperbolic.
The main result of this section is the following theorem: 
\begin{thm}\label{mainthm}
For $s>3/2$ and under the assumptions \rm{(A1), (A2)} and \rm{(A3)}, if $W_0 \in H^s$ satisfy the non-emptiness condition
\begin{equation}
\inf_{x \in \mathds{R}} {h}_0(x)\ =\ \inf_{x \in \mathds{R}} \left( \eta_0(x)\, +\, 
d(0,x) \right)\ \geqslant\ h^* >\ 0,
\end{equation}
then there exist $T>0$ and a unique solution $W \in {C}([0,T], H^s) \cap {C}^1([0,T], 
H^{s-1})$ of the system \eqref{generalsystem}. Moreover, if the maximal existence time $T_{max}< +\infty$,  
then 
\begin{equation}\label{buc}
\inf_{(t,x) \in [0,T_{max}[ \times \mathds{R}} h(t,x)\ =\ 0 \qquad \text{or} \qquad \lim_{t \to T_{max}} 
\|W\|_{H^s}\ =\ +\infty.
\end{equation}
\end{thm}

{\bf Remarks:}
(i) Theorem \ref{mainthm} holds also for periodic domains; 
(ii) The right-hand side of \eqref{generalsystem} can be replaced by a finite sum on the form
\begin{equation}\label{sumop}
F(W)\ =\  \begin{pmatrix} \mathfrak{A}_1\, f_1\ +\  \mathfrak{A}_3\, f_3 \\  \mathfrak{A}_2\, f_2\ +\  \mathfrak{A}_4\, f_4 \end{pmatrix}\ +\ \begin{pmatrix} \mathfrak{B}_1\, k_1\ +\  \mathfrak{B}_3\, k_3 \\  \mathfrak{B}_2\, k_2\ +\  \mathfrak{B}_4\, k_4 \end{pmatrix}\ +\ \cdots,
\end{equation}
where the additional terms satisfy also the conditions (A1), (A2) and (A3); 
(iii) Under some additional assumptions, the blow-up criteria \eqref{buc} can be improved (see 
Theorem \ref{thm:blowup}, for example); 
(iv) If for some $2 \leqslant i \leqslant N$, the function $g_i$ appears only on $f_1$ and $f_2$, 
then, due to (A2), the assumption $g_i \in {C}(\mathds{R}^+,H^s)$ can be replaced by $g_i \in 
{C}(\mathds{R}^+,H^{s-1})$. \\

In order to prove the local well-posedness of \eqref{generalsystem}, we consider 
\begin{equation}\label{iteration}
\partial_t\, W^{n+1}\ +\ B(W^n)\, \partial_x\, W^{n+1}\ =\ F(W^n), \qquad 
W^n(0,x)\ =\ (\eta_0(x),u_0(x))^T,
\end{equation}
where $n\geqslant0$ and $W^0(t,x)=(\eta_0(x),u_0(x))^\top$. 
The idea of the proof is to solve the linear system \eqref{iteration}, then, taking the limit $n \to 
\infty$, we obtain a solution of \eqref{generalsystem}.
Note that we have assumed that  $g_1$ and $f$ are positive, so $g_1 f > 0$, then the system 
\eqref{iteration} is hyperbolic; it is an important point to solve each iteration in \eqref{iteration}.
Note that a symmetriser of the matrix $B(W)$ is 
$A(W)\ \eqdef\ \begin{pmatrix} g_1\, f(d,h) & 0 \\ 0 & 1 \end{pmatrix}$.

Let $(\cdot,\, \cdot)$ be the scalar product in $L^2$ and let the energy of \eqref{iteration} be defined as 
\begin{equation*}
%
E^{n+1}(t)\ \eqdef\ \left(\Lambda^s\, W^{n+1},\ A^n\, \Lambda^s\, W^{n+1} \right)\ \forall t \geqslant 0.
\end{equation*}
If $g_1 f$ is bounded and far from $0$, then $E^n(t)$ is equivalent to $\|W^n\|_{H^s}$. In order 
to prove Theorem \ref{mainthm}, the following results are needed.
\begin{thm}\label{energyestimate}
Let $s>3/2$, ${h}^{*}>0$ and $R>0$, then there exist $K,T>0$ such that: if the initial data 
$(\eta_0,{h}_0) \in H^s$ satisfy
\begin{equation}\label{initial}
\inf_{x \in \mathds{R}} {h}_0(x)\ \geqslant\ 2\, {h}^{*}, \qquad E^n(0)\ <\ R,
\end{equation}
and $W^n \in {C}([0,T], H^s) \cap {C}^1([0,T], H^{s-1})$, satisfying for all $t \in [0,T]$
\begin{equation}
h^n\ \geqslant\ {h}^{*}, \qquad \|(W^n)_t\|_{H^{s-1}}\ \leqslant\ K, \qquad E^n(t)\ \leqslant\ R,
\end{equation}
then there exists a unique $W^{n+1} \in C([0,T], H^s) \cap C^1([0,T], H^{s-1})$ solution of 
\eqref{iteration} such that 
\begin{equation}
h^{n+1}\ \geqslant\ h^{*}, \qquad \|(W^{n+1})_t\|_{H^{s-1}}\ \leqslant\ K, \qquad E^{n+1}(t)\ \leqslant\ R.
\end{equation}
\end{thm}
The proof of Theorem \ref{energyestimate} is classic (it can be done following \citet{guelmame2020Euler,Israwi2011,liu2019well} and 
using the following lemmas).

Let $\Lambda$ be defined such that $\widehat{\Lambda f}=(1+\xi^2)^\frac{1}{2} \hat{f}$, and let  
$[A,B] \eqdef AB-BA$ be the commutator of the operators $A$ and $B$. We have the following lemma.
\begin{lem}(\citet{kato1988commutator})
If $r\, \geqslant\, 0$, then 
\begin{align}
\|f\, g\|_{H^r}\ &\lesssim\ \|f\|_{L^\infty}\, \|g\|_{H^r}\ +\ \|f\|_{H^r}\, \|g\|_{L^\infty}, 
\label{Algebra} \\
\left\| \left[ \Lambda^r,\, f \right]\, g  \right\|_{L^2}\ &\lesssim\  \|f_x\|_{L^\infty}\, \|g\|_{H^{r-1}}\ 
+\ \|f\|_{H^r}\, \|g\|_{L^\infty}. \label{Commutator}
\end{align}
\end{lem}

\begin{lem}
Let $k \in \mathds{N}^*$, $F \in {C}^{m+2}(\mathds{R}^k)$ with $F(0,\cdots,0)=0$ and $0 \leqslant s 
\leqslant m$, then there exists a continuous function $\tilde{F}$, such that for all $f=(f_1,\cdots,f_k) 
\in H^s \cap W^{1,\infty}$ we have
\begin{equation}\label{Composition2}
\|F(f)\|_{H^s}\ \leqslant\ \tilde{F} \left( \|f\|_{W^{1,\infty}} \right)\, \|f\|_{H^s}.
\end{equation}
\end{lem}
\proof
The case $k=1$ has been proved in \cite{constantin2002initial}. Here, we prove the inequality \eqref{Composition2} by induction (on $s$). Note that {\small
\begin{gather*}
F(f_1,\cdots,f_k)\,
=\, F(0,f_2,\cdots,f_k)\ +\ \int_0^{f_1}F_{f_1}(g_1,f_2,\cdots,f_k)\, \ud g_1  \nonumber\\
=\, F(0,0,f_3,\cdots,f_k)\, +\, \int_0^{f_1}\!  F_{f_1}(g_1,f_2,\cdots,f_k)\, \ud g_1\, +\, 
\int_0^{f_2}\! F_{f_2}(0,g_2,f_3,\cdots,f_k)\, \ud g_2\ +\ \cdots \nonumber\\
=\, \int_0^{f_1}F_{f_1}(g_1,f_2,\cdots,f_k)\, \ud g_1\, +\, \cdots\, +\, \int_0^{f_k}F_{f_k}(0,
\cdots,0,g_k)\, \ud g_k .
\end{gather*}}
This implies that 
\begin{equation}\label{L2}
\|F(f_1,\cdots,f_k)\|_{L^2}\ \lesssim\ \|f\|_{L^2},
\end{equation}
which is \eqref{Composition2} for $s=0$. For $s \in ]0,1[$, let 
\begin{align*}
\big| F(f_1&(x+y),\cdots,f_k(x+y))\, -\, F(f_1(x),\cdots,f_k(x)) \big|\\ 
\leqslant &\ \big| F(f_1(x+y),\cdots,f_k(x+y))\, -\, F(f_1(x),f_2(x+y),\cdots,f_k(x+y)) \big|\\
&+\ \big| F(f_1(x),f_2(x+y),\cdots,f_k(x+y))\, -\, F(f_1(x),f_2(x),f_3(x+y),\cdots,f_k(x+y)) \big|\\
&+\ \cdots\ +\ \big| F(f_1(x),\cdots,f_{k-1}(x),f_k(x+y))  -\, F(f_1(x),\cdots,f_k(x))\big|\\
\leqslant &\ \sum_{i=1}^k\, \left| f_i(x+y)\, -\, f_i(x) \right|\, \|F_{f_i}\|_{L^\infty}. 
\end{align*}
The last inequality, with the definition $H^s \eqdef \left\{f \in L^2,\ \int_\mathds{R}\int_\mathds{R} \frac{|f(x+y)-f(x)|^2}{|y|^{1+2\, s}}\, \ud x \,\ud y < +\infty \right\}$, 
implies \eqref{Composition2} for $s \in ]0,1[$. 
For $s \geqslant 1$, the proof is done by induction. 
Using \eqref{Algebra} and \eqref{L2}, we obtain
$\|F(f)\|_{H^s}\/ 
\lesssim\/ \left\| \sum_{i=1}^k\, F_{f_i}(f)\, \partial_x\, f_i \right\|_{H^{s-1}}\/ +\/ 
\|F(f)\|_{L^2}\/\lesssim\/ \|f\|_{H^s}\/ +\/ \sum_{i=1}^k\, \| F_{f_i}(f)\|_{H^{s-1}}$.
Using the induction and the last inequality, we obtain \eqref{Composition2} for all $s \geqslant 0$. 
\qed \\


\textbf{Proof of Theorem \ref{mainthm}.}
Using Theorem \ref{energyestimate}, one obtains that $(W^n)$ is uniformly bounded in $ {C}([0,T], H^s) 
\cap {C}^1([0,T], H^{s-1})$ and satisfies ${h}^n \geqslant {h}^{*}$. Defining
\begin{equation}
\tilde{E}^{n+1}(t)\ \eqdef\ \left( \Lambda^{s-1}\, (W^{n+1}\, -\, W^n),\ A^n\, \Lambda^{s-1}\, 
(W^{n+1}\, -\, W^n) \right),
\end{equation}
and using \eqref{iteration}, one obtains
\begin{align}
\tilde{E}_t^{n+1}\ 
%
=&\ 2\, \left( \Lambda^{s-1} \left( F^n-F^{n-1} + (B^{n-1}-B^n)\/W_x^n \right) ,\ A^n\, 
\Lambda^{s-1}\, (W^{n+1}\, -\, W^n) \right) \nonumber\\
&-\ 2\, \left( [\Lambda^{s-1},B^n] (W^{n+1}-W^n)_x,\ A^n\, \Lambda^{s-1} (W^{n+1}-W^n) \right) 
\nonumber\\
&+\ \left( \Lambda^{s-1} \left( (A^n B^n)_x (W^{n+1}-W^{n}) \right) ,\ \Lambda^{s-1}\, (W^{n+1}\, 
-\, W^n) \right)     \nonumber\\
&+\ \left( \Lambda^{s-1}\, (W^{n+1}\, -\, W^n),\ (A^n)_t\, \Lambda^{s-1}\, (W^{n+1}\, -\, W^n)  \right).
\end{align}
Using (A2), (A3), \eqref{Algebra} and \eqref{Composition2}, one obtains
\begin{equation}
\|F^n-F^{n-1}\|_{H^{s-1}}\ +\ \|(B^{n-1}-B^n)\partial_x W^n\|_{H^{s-1}}\ \lesssim\ \|W^n - W^{n-1}\|_{H^{s-1}}\ \lesssim\ \sqrt{\tilde{E}^{n}},
\end{equation}
where ``$\mathscr{A}\lesssim \mathscr{B}$'' means $\mathscr{A} \leqslant C \mathscr{B}$, with $C>0$ is a constant independent of $n$. 
Using \eqref{Commutator} and \eqref{Composition2}, we obtain
\begin{equation}
\left\|[\Lambda^{s-1},B^n] (W^{n+1}-W^n)_x \right\|_{L^2}\ \lesssim\ \|W^{n+1} - W^{n}\|_{H^{s-1}}\  
\lesssim\ \sqrt{\tilde{E}^{n+1}}.
\end{equation}
From \eqref{Algebra} and \eqref{Composition2}, it follows that 
\begin{equation}
\|(A^n B^n)_x (W^{n+1}-W^{n})\|_{H^{s-1}}\ \lesssim\ \|W^{n+1} - W^{n}\|_{H^{s-1}}\  \lesssim\ 
\sqrt{\tilde{E}^{n+1}}.
\end{equation}
Combining the estimates above, we obtain that
$\tilde{E}_t^{n+1}\lesssim\tilde{E}^{n+1}+ \tilde{E}^{n}$,
and using that $\tilde{E}^n(0)=0$, we obtain 
$\tilde{E}^{n+1} \leqslant \left( \mathrm{e}^{C\/t}\, -\, 1 \right) \tilde{E}^{n}$ where $C>0$ does not depend on $n$.
Taking $T>0$ small enough, it follows that 
\begin{equation}
\|W^{n+1}-W^n\|_{H^{s-1}} \lesssim\ \tilde{E}^{n+1}\ \leqslant\ \half\, \tilde{E}^n\ \leqslant\ {\textstyle \frac{1}{2^{n}}}\, \tilde{E}^1. 
\end{equation}
Finally, taking the limit $n \to \infty$ in the weak formulation of \eqref{iteration} and using (A3), 
we obtain a solutions of the system \eqref{system}. This completes the proof of Theorem \ref{mainthm}.
\qed

\section{Proof of Theorems \ref{existenceSV} and \ref{thm:blowup}}\label{sec:proof}

The system \eqref{rSVmb2} is written in the form \eqref{system} by replacing the right-hand side 
of \eqref{system}, as in \eqref{sumop}, taking $N = 4$ and 
$G(t,x)=(g,d_x,d_{xx},d_t)$, $a(d,h,u,g,d_x,d_{xx},d_t)=c(d,h,u,g,d_x,d_{xx},d_t)=u$,  
$b(d,h,u,g,d_x,d_{xx},d_t)=h$, $f(d,h)=h^{-1}$, $f_1=f_4=k_1=k_3=k_4=0$,  
$f_2(d,h,u,h_x,u_x,g,$ $d_x,d_{xx},d_t)=2\/h^3\/u_x^{\,2}\/-\/\half\/g\/h^2(\eta_x^{\,2}+
2\/\eta_x\/d_x)$, $f_3(d,h,u,g,d_x,d_{xx},d_t)=-d_t\/-u\/d_x$, $k_2(d,h,u,h_x,u_x,$ $g,d_x,
d_{xx},d_t)= g\/h^2\/\eta_x\/d_{xx}$, $\mathfrak{A}_1=\mathfrak{A}_4=\mathfrak{B}_1=\mathfrak{B}_3=\mathfrak{B}_4=0$, $\mathfrak{A}_2=- \varepsilon\/\mathcal{L}_h^{-1}\/
\partial_x$ and $\mathfrak{A}_3=1$, $\mathfrak{B}_2= \varepsilon\/\mathcal{L}_h^{-1}$.
Then, in order to prove Theorem \ref{existenceSV}, the following lemma is needed:
\begin{lem}(\citet{liu2019well})\label{Inverseesitimates}
Let $0< h_{\mathrm{inf}} \leqslant h \in W^{1,\infty}$, then the operator $\mathcal{L}_h$ is an isomorphism from $H^2$ to $L^2$ and 
if $0 \leqslant s \leqslant \tilde{m} \in \mathds{N}$, then 
\begin{equation}\label{estimate2}
\left\|\mathcal{L}_h^{-1}\, \psi \right\|_{H^{s+1}}\ +\
\left\|\mathcal{L}_h^{-1}\, \partial_x\, \psi \right\|_{H^{s+1}}\  \leqslant\ C\, \left\|\psi \right\|_{H^s}\, \left(1 +\ \left\|h\, -\, \bar{d} \right\|_{H^s} \right), 
\end{equation}
where $C$ depends on $s$, $\varepsilon$, $h_\mathrm{inf}$, $\|h - \bar{d}\|_{W^{1,\infty}}$ and not on 
$\|h - \bar{d}\|_{H^s}$.
\end{lem}

\textbf{Proof of Theorem \ref{existenceSV}.}
In order to prove Theorem \ref{existenceSV}, it suffices to verify (A1)--(A3). The assumption (A1) 
is obviously satisfied and (A2) follows from Lemma \ref{Inverseesitimates}.
In order to prove (A3), let $W,\tilde{W},\psi \in H^s$. Using Lemma \ref{Inverseesitimates} and 
\eqref{Algebra}, we obtain
\begin{gather*}
\left\| \left( \mathcal{L}_h^{-1}\, -\, \mathcal{L}_{\tilde{h}}^{-1} \right) \psi \right\|_{H^{s-1}}\ 
=\,   \left\| \mathcal{L}_h^{-1}\, \left(  \mathcal{L}_{\tilde{h}}\, -\, \mathcal{L}_h\, \right) 
\mathcal{L}_{\tilde{h}}^{-1} \psi \right\|_{H^{s-1}}\\
\lesssim\,\left\| \left(  \mathcal{L}_{\tilde{h}}\, -\, \mathcal{L}_h\, \right) \mathcal{L}_{\tilde{h}}^{-1} \psi \right\|_{H^{s-2}}\ 
\lesssim\,  \left\| h\, -\, \tilde{h} \right\|_{H^{s-1}}\ \leqslant\, \left\|W\, -\, \tilde{W} 
\right\|_{H^{s-1}}.
\end{gather*}
where the constants depend on $s,d,\|W\|_{H^s},\|\tilde{W}\|_{H^s},\|\psi\|_{H^{s-1}} $. 
The same proof can be used with the operator $\mathfrak{A}_2$. \qed

\textbf{Proof of Theorem \ref{thm:blowup}.} Using the characteristics $\chi(0,x)=x$ and 
$\chi_t(t,x)=u(t,\chi(t,x))$, the conservation of the mass \eqref{mc} becomes 
\begin{equation}\label{h}
{\ud h}/{\ud\/t}\ +\ u_x\, {h}\  =\ 0, \qquad \implies \qquad   {h}_0(x)\, \mathrm{e}^{-t\, 
\|u_x\|_{L^\infty}}\ \leqslant\ {h}(t,x)\ \leqslant\ {h}_0(x)\, \mathrm{e}^{t\, \|u_x\|_{L^\infty}}.
\end{equation} 
The energy equation \eqref{eneSV} implies that 
\begin{equation} 
\mathscr{E}'(t)\ \leqslant\ \left( |\dot{g}|/g\, +\, 1 \right)\, \mathscr{E}(t)\ +\ \half\, g\/ 
\int_\mathds{R} \left(\,d_t^{\,2}\, +\, \varepsilon\, h^2\, d_{xt}^{\,2}\, \right)  \ud\/x,
\end{equation}
since $h$ is bounded, the inequality \eqref{eb} follows by Gronwall's lemma.

In order to prove the blow-up criterion, we first suppose that $\|W_x\|_{L^\infty}$ is bounded and 
we show that the scenario \eqref{buc1} is impossible.
The equation \eqref{h} implies that $h$ is bounded and far from $0$. Using $\|W\|_{L^\infty} \leqslant 
\|W\|_{H^1} \lesssim \mathscr{E}(t)$, one obtains that $\|W\|_{W^{1,\infty}}$ is bounded on any 
interval $[0,T]$.
Using Lemma \ref{Inverseesitimates} and doing some classical energy estimates (see \cite{guelmame2020Euler,Israwi2011,liu2019well}), 
we can prove that $\|W\|_{H^s}$ is also bounded. This ends the proof of Theorem \ref{thm:blowup}.  \qed

%



\end{document}